\documentclass{jhrs}
\usepackage{amsmath}
\usepackage{xy}

\begin{document}
\title{Rational formality of mapping spaces}
\author{Yves F\'elix}
\email{yves.felix@uclouvain.be}
\address{Universit\'e Catholique de Louvain\\
Institut de Math\'ematique et de Physique\\
2, chemin du cyclotron\\
1348 Louvain-La-Neuve\\ Belgium}


\begin{abstract}
Let $X$ and $Y$ be finite nilpotent   CW complexes with dimension of
$X$  less than the connectivity of $Y$. Generalizing results of
Vigu\'e-Poirrier and Yamaguchi, we prove that the mapping space
$\mbox{Map}(X,Y)$ is rationally formal if and only
 if $Y$ has the rational homotopy type of a finite product of odd dimensional spheres.
\end{abstract}

\received{September 29, 2009}   
\published{Month Day, Year}  
\submitted{Pascal Lambrechts}  

\volumeyear{2006} 
\volumenumber{1}  
\issuenumber{1}   

\startpage{1}     
\maketitle

\section{Introduction}
Let $X$ and $Y$ be connected spaces that have the rational homotopy
type of finite CW complexes. We denote by $n$ the maximum integer
$q$ such that $H^q(X;\mathbb Q)\neq 0$. In this text we consider
mapping spaces $\mbox{Map}(X,Y)$ satisfying the following
hypotheses~(H).
$$\leqno{H}\left\{ \begin{array}{l}
  \mbox{ $(i)$   $X$ and $Y$ are not rationally contractible},\\
\mbox{ $(ii)$  There exists $n\geq 1$ such that $H^n(X;\mathbb Q)\neq 0$, $H^q(X;\mathbb Q)= 0$ if $q>n$}, \\\mbox{} \hspace{1cm} \mbox{  and    $Y$ is $n$-connected }
     \end{array}\right.$$

     \vspace{2mm}\noindent Under those hypotheses, $\mbox{Map}(X,Y)$ is a nilpotent space
     and its rational homotopy is described by Haefliger \cite{Ha} and Brown and
     Szczarba \cite{BS}.

Our main interest here is to understand when $\mbox{Map}(X,Y)$ is a
(rationally) formal space.
  Formality is   important in rational homotopy. If a
       space is formal then its rational homotopy type is completely determined
       by its rational cohomology. More precisely a nilpotent space $Z$ is formal
        if its Sullivan minimal model is quasi-isomorphic to  the
         differential graded algebra $(H^*(Z;\mathbb Q),0)$. Many spaces coming from geometry are formal.
         Among formal spaces we find
         the spheres, the projective spaces, the products of Eilenberg-MacLane
         spaces, the compact K\"ahler manifolds   (\cite{DGMS}),
           and the $(p-1)$-connected compact manifolds, $p\geq 2$, of dimension $\leq 4p-2$
       \cite{Mi}.

    The formality of mapping spaces has been the subject of previous works. In \cite{DV},
    N. Dupont and M. Vigu\'e-Poirrier prove that when $H^*(Y;\mathbb Q)$ is finitely generated, then
    $\mbox{Map}(S^1,Y)$  is formal if and only if $Y$ is
    rationally a product of Eilenberg-MacLane spaces. In \cite{Ya} T. Yamaguchi proves
     that when $Y$ is elliptic, the formality of $\mbox{Map}(X,Y)$ implies that $Y$
     is rationally a product of odd dimensional spheres. In \cite{Vig} M. Vigu\'e-Poirrier
     proves that if $\mbox{Map}(X,Y)$ is formal and if the Hurewicz map
     $\pi_q(X)\otimes \mathbb Q \to H_q(X;\mathbb Q)$ is nonzero in some odd degree $q$,
     then $Y$ has the homotopy type of a product of Eilenberg-MacLane spaces.
     When $Y$ is a finite complex, we prove here that   the hypothesis  on the Hurewicz map is not necessary.

     \vspace{2mm}\noindent {\bf Theorem 1.} {\sl Under the above hypotheses (H),
     $\mbox{Map}(X,Y)$ is formal if and only if $Y$ has the rational homotopy type
     of a product of odd dimensional spheres.}

     \vspace{2mm}As an important step in the proof of Theorem 1 we prove

     \vspace{2mm}\noindent {\bf Theorem 2.} {\sl If dim$\, Y=N$, then the Hurewicz map
     $$\pi_q(\mbox{Map}(X,Y))\otimes \mathbb Q \to H_q(\mbox{Map}(X,Y);\mathbb Q)$$
     is zero for $q>N$. }

     \section{Rational homotopy}

     The theory of minimal models originates in the works of Sullivan \cite{Su} and Quillen \cite{Qu}.
  For recall
     a graded algebra $A$ is graded commutative if $ab= (-1)^{\vert a\vert \cdot \vert a\vert} ba$
     for homogeneous elements $a$ and $b$. A graded commutative algebra $A$ is free on a graded
     vector space $V$, $A = \land V$, if $A$ is the quotient of the tensor algebra $TV$
     by the ideal generated by the elements $xy- (-1)^{\vert x\vert \cdot \vert y\vert} yx$,
     $x,y\in V$. A  (Sullivan) minimal algebra is a graded commutative differential   algebra
     of the form $(\land V,d)$ where $V$ admits a basis $v_i$ indexed by a well ordered
      set $I$ with $d(v_i)\in \land (v_j, j<i)$. Now if $(A,d)$ is a graded
      commutative differential   algebra whose cohomology is
      connected and finite type, there is a unique (up to
      isomorphism) minimal   algebra $(\land V,d)$ with a
      quasi-isomorphism $\varphi : (\land V,d) \to (A,d)$. The
      differential graded algebra $(\land V,d)$ is then called the
      (Sullivan) minimal model of $(A,d)$.

     In \cite{Su} Sullivan associated to each nilpotent space $Z$ a graded commutative differential   algebra
      of rational
     polynomials forms on $Z$, $A_{PL}(Z)$, that is a rational replacement of the algebra
     of de Rham forms on a manifold. The  minimal model $(\land V,d)$
     of $ A_{PL}(Z)$ is then called the minimal model of $Z$.   More generally a model
      of $Z$ is a  graded commutative differential   algebra quasi-isomorphic to its
       minimal model.   For more details we refer to \cite{Su}, \cite{FHT} and \cite{FOT}.

     A space $X$ is called (rationally) formal if its minimal model, $(\land V,d)$, is quasi-isomorphic to
     its cohomology with differential $0$,
     $$\psi : (\land V,d) \to (H^*(X;\mathbb Q), 0)\,.$$
    A formal space $X$ admits a minimal model equipped with a  bigradation on $V$, $V =
    \oplus_{p\geq 0,q\geq 1} V_p^q$ such that $d(V_p^q)\subset  (\land  V)_{p-1}^{q+1}$,
    and such that  the bigradation induced on the homology satisfies $H_p^q= 0$ for $p\neq 0$.
    This model has been constructed by Halperin and Stasheff in \cite{HS}, and is called the
     bigraded model of $X$. We will use this model for the proof of Theorem 2.

     A nilpotent space $X$ is called (rationally) elliptic if $\pi_*(X)\otimes \mathbb Q$ and
      $H^*(X;\mathbb Q)$ are finite dimensional vector spaces. To be elliptic for a space
      $X$ is a very restrictive condition. For instance $H^*(X;\mathbb Q)$ satisfies Poincar\'e
      duality and $\pi_q(X)\otimes \mathbb Q$ is zero for $q\geq$ $2\cdot$ dim$\, X$. A
      nilpotent space $X$ is called (rationally) hyperbolic if $\pi_*(X)\otimes \mathbb Q$
      is infinite dimensional and $H^*(X;\mathbb Q)$ finite dimensional. The homotopy groups of
      elliptic and hyperbolic spaces have a completely different behavior. For instance,
      for an hyperbolic space $X$,  the
      sequence $\sum_{i\leq q} \dim\, \pi_i(X)\otimes \mathbb Q$ has an exponential growth (\cite{FHT}).

     In \cite{Ha}, Haefliger gives a process to construct a minimal model for $\mbox{Map}(X,Y)$.
     With the hypotheses (H) of the Introduction, suppose that $(\land W,d)$ is the Sullivan
     minimal model of $X$. Denote by  $S\subset (\land W)^n$
    a supplement of the subvector space generated by the cocycles. Then $I = (\land W)^{>n} \oplus S$
     is an acyclic
     differential graded ideal, and the quotient $(A,d)= (\land W/I,d)$ is a finite dimensional
     model for $X$. We denote by $(B,d)$ the dual coalgebra. Let
     $(a_i)$, $i=0, \ldots ,p$ be a graded basis for $A$ with $a_0=1$ and denote by $\overline{a_i}$ the
     dual basis for $B$.

 Denote also by $(\land V,d)$   the minimal model of $Y$. We define a morphism of graded algebras
     $$\varphi: \land V \to A\otimes \land (B\otimes V)$$
     by putting
     $\varphi (v) = \sum_i a_i \otimes (\overline{a_i}\otimes v)$.
     In \cite{Ha} Haefliger proves that there is a unique differential $D$ on $\land (B\otimes V)$ making
     $$\varphi : (\land V,d)\to (A,d)\otimes (\land (B\otimes V),D)$$
     a morphism of differential graded algebras. Then $(\land (B\otimes V),D)$ is a model for
     $\mbox{Map}(X,Y)$ and $\varphi$ is a model for the evaluation map $\mbox{Map}(X,Y)\times X\to Y$.
     In particular, (\cite{Vi}), the rational homotopy groups of $\mbox{Map}(X,Y)$ are given by
     $$\pi_q(\mbox{Map}(X,Y))\otimes\mathbb Q = \oplus_i \left[H_i(X;\mathbb Q)\otimes
     \pi_{q+i}(Y)\otimes \mathbb Q\right]\,.$$ This formula is natural in $X$ and   $Y$.

\section{Proof of Theorem 1.} In \cite{Th} Thom computes the rational homotopy type of
$\mbox{Map}(X, K(\mathbb Q,r))$ when dim$\, X<r$. He proves that the mapping space is a product of
 Eilenberg-MacLane spaces,
$$\mbox{Map}(X,K(\mathbb Q,r)) = \prod_i K(H_i(X;\mathbb Q), r-i)\,.$$
Since odd dimensional spheres are rationally Eilenberg-MacLane spaces, it follows that if $Y$ has
the rational homotopy type of a product of odd dimensional spheres, then $\mbox{Map}(X,Y)$ is formal.

Suppose now that $\mbox{Map}(X,Y)$ is formal. Since any retract of a formal space is formal, $Y$ is
 formal. By Theorem 2, the image of the Hurewicz map for $\mbox{Map}(X,Y)$ is finite dimensional.
 Recall that  for a formal space, the
 cohomology is generated by classes that evaluate non trivially on the image of the Hurewicz map.
 Therefore the algebra $H^*(\mbox{Map}(X,Y);\mathbb Q)$ is finitely generated.

The square of an even dimensional generator $x_i$ of
$H^*(\mbox{Map}(X,Y);\mathbb Q)$  gives a map $\mbox{Map}(X,Y) \to
K(\mathbb Q, 2r_i)$, $r_i = 2\vert x_i\vert$. We denote by $\theta$
the product of those maps,
 $$\theta : \mbox{Map}(X,Y)\to \prod_i K(\mathbb Q, 2r_i)\,.$$ We do not suppose that $x_i^2\neq 0$. In fact if $x_i^2= 0$ for all $i$, then $\theta$ is homotopically trivial but this has no effect on our argument.
 The pullback along $\theta$ of the product of the principal fibrations $K(\mathbb Q, 2r_i-1)\to
  PK(\mathbb Q, 2r_i)\to K(\mathbb Q, 2r_i)$ is a fibration
 $$\prod_i \, K(\mathbb Q, 2r_i-1) \to E \to \mbox{Map}(X,Y)\,.$$
 By construction the rational cohomology of $E$ is finite dimensional, and so the rational category
  of $E$ is also finite.

 Now from the definition of the dimension of $X$, there is a cofibration
 $X' \to X \stackrel{q}{\to} S^n$ such that $H_n(q;\mathbb Q)$ is surjective. The restriction to
 $X'$ induces a map $\mbox{Map}(X,Y)\to \mbox{Map}(X',Y)$ whose homotopy fiber is the injection
 $$j : \Omega^n Y= \mbox{Map}_*(S^n,Y)\to \mbox{Map}(X,Y)\,.$$ From the naturality of the formula
  for the rational homotopy groups of a mapping space, we deduce that $\pi_*(j)\otimes \mathbb Q$
  is injective.
 Denote now $E'$ the pullback of $E\to \mbox{Map}(X,Y)$ along $j$,
 $$\begin{array}{ccc}
 \prod_iK(\mathbb Q, 2r_i-1) & =& \prod_i \, K(\mathbb Q, 2r_i-1)\\
 \downarrow &&\downarrow \\
 E' & \stackrel{j'}{\to} & E\\
 \downarrow &&\downarrow \\
 \Omega^n Y & \stackrel{j}{\to} & \mbox{Map}(X,Y)\end{array}$$
Since $\pi_*(j')\otimes \mathbb Q$ is injective, it   follows from
the mapping theorem  \cite{FHT}  that the rational category of $E'$
is finite. In particular the cup length of $E'$ is finite.

Now the rational cohomology of $\Omega^n Y$ is the free commutative graded algebra on the graded vector
 space $S_*$, with $S_q = \pi_{n+q}(Y)\otimes \mathbb Q$. Therefore if $Y$ is hyperbolic,
 $H^*(E';\mathbb Q)$ will contain a free commutative graded algebra on an infinite number of
  generators, and in particular its cup length is infinite. It follows that $Y$ is elliptic.
  To end the proof we only apply Yamaguchi result (\cite{Ya}) that asserts that when $Y$ is elliptic,
  and $\mbox{Map}(X,Y)$ is formal, then $Y$ has the rational homotopy type of a finite product of odd
   dimensional spheres.

\section{Proof of Theorem 2}
Denote by $(\land V,d)$ the bigraded model for $Y$ and by $(A,d)$ a
connected finite dimensional model for $X$. Connected means that
$A^0 = \mathbb Q$. Denote as above by $a_i$, an homogeneous basis of
$A$, and by $\overline{a_i}$ the dual basis for $B =
\mbox{Hom}(A,\mathbb Q)$. We write also $B_+ =
\mbox{Hom}(A^+,\mathbb Q)$.

Recall now that a model for the evaluation map $X\times \mbox{Map}(X,Y)\to Y$ is given by the morphism
$$\varphi : (\land V,d) \to (A,d)\otimes (\land (B\otimes V),D)\,,$$
defined by $\varphi (v) = \sum_i a_i\otimes (\overline{a_i}\otimes v)$.

We consider the differential ideal  $I= \land V \otimes \land^{\geq
2} (B_+\otimes V) $, and we denote by
 $\pi : (\land (B\otimes V),D)\to (\land (B\otimes V)/I, \overline{D})$ the quotient map.  In
 $\land (B\otimes V)/I$ the equation $\pi\circ \varphi\circ d=
 (d\otimes 1+ 1\otimes \overline{D})\circ \pi\circ \varphi$ gives for
 each $v\in V$ the equation
$$\sum_i \, da_i \otimes (\overline{a_i}\otimes v) + \sum_i (-1)^{\vert a_i\vert} a_i\otimes
 \overline{D} \,(\overline{a_i}\otimes v) = 1\otimes dv + \sum_{a_i\in A^+} a_i \otimes \theta_i(v)\,,$$
where $\theta_i$ is the derivation of $\land V \otimes \land (B\otimes V)$ defined by $\theta_i(v) =
\overline{a_i}\otimes v$ and $\theta_i (B\otimes V) = 0$.

To go further we specialize the basis of $A^+$. We denote by
$\{y_i\}$ a basis of $d(A)$, by $\{ e_j\}$ a set of cocycles such
that $\{y_i, e_j\}$ is a basis of the cocycles in $A$. Finally we
choose elements $x_i$ with $d(x_i)= y_i$. A basis of $A$ is then
given by $1$ and  the elements $ x_i, y_i$ and $e_j$. Denote then by
$\psi_j$, $\psi'_i$ and $\psi''_i$ the derivations $\theta $
associated respectively to $e_j$, $x_i$ and $y_i$. Then we have
$$\overline{D} (\overline{e_i}\otimes v) = (-1)^{\vert e_i\vert} \psi_i (v)\,, \hspace{3mm}
\overline{D}(\overline{x_i}\otimes v) = (-1)^{\vert x_i\vert} \psi'_i(v)\,,$$
$$\overline{D}(\overline{y_i}\otimes v) = (-1)^{\vert y_i\vert} \left( \psi''_i(v) - (\overline{x_i}
\otimes v)\right)\,.$$

it follows  that the complex $(\land (B\otimes V)/I, \overline{D})$
decomposes into a direct sum
$$(\land (B\otimes V)/I, \overline{D}) = \land V \oplus  ( \oplus_j C_j)\oplus D \,,
\hspace{1cm}\mbox{with }  C_j =( \overline{e_j}\otimes V)\otimes
\land V\,,$$ and where $D$ is the ideal generated by the
$\overline{x_i}\otimes v$ and
 $\overline{y_i}\otimes v$.

Consider now in $(\land (B\otimes V), D)$ a cocycle $\alpha$ of the
form $$\alpha = \sum_j \overline{e_j} \otimes v_j + \sum_i
\overline{x_i}\otimes u_i + \sum_i \overline{y_i}\otimes w_i+
\omega$$ where $\omega$ is a decomposable element. Looking at the
linear term of $D(\alpha)$ we obtain that $\sum_i \overline{y_i}
\otimes w_i= 0$. We can replace $\alpha$ by $\alpha + D(\sum_i
(-1)^{\vert x_i\vert}\overline{y_i} \otimes u_i)$ to cancel the
linear part $\sum_i \overline{x_i}\otimes u_i$. We can thus suppose
that $\alpha$ has the form
$$\alpha = \sum_j \overline{e_j}\otimes v_j + \omega$$ where $\omega$
is a decomposable element.

In $\land (B\otimes V)/I$,  $\alpha$ decomposes into a sum of
cocycles, $\alpha =
 \sum_i \alpha_i$
 with $\alpha_i\in C_i$.
Let fix some $i$. We write $r = \vert e_i\vert$ and  $\overline{v}=
(\overline{e_i}\otimes v)$. We denote $\overline{V}=
\overline{e_i}\otimes V$. Then the component $C_i$ is isomorphic to
$ (\land V \otimes
  \overline{V}, \overline{D})\, $ and $\overline{V}$ is   equipped with an
isomorphism of degree $-r$,
$$s : V^q\to \overline{V}^{q-r}\,.$$
We extend $s$ in a derivation of $\land V\otimes \land \overline{V}$
by $s(\overline{V})=0$, and the differential $\overline{D}$
satisfies $\overline{D}(\overline{v}) = (-1)^r sd(v)$.

Write $\alpha_i = \overline{v} + \omega$, where $\omega  \in
\overline{V}\otimes \land^+V$. We show that in that case $v$ is a
cocycle. If this is true for any $i$, this implies that the map
$$\rho_q: H^q(\land V\otimes \land (B\otimes V),D) \to H^q((\land
 V\otimes \land (B\otimes V) ) /\land^{\geq 2} (V\oplus
(B\otimes V))),D) $$ is zero in degrees $q\geq$ dim$\, Y$. Since
$\rho_q$ is the dual of the Hurewicz map $h_q : \pi_q(\mbox{Map}(X,
Y))\otimes \mathbb Q \to H_q(\mbox{Map}(X, Y);\mathbb Q)$, this
implies the result.

 We now
follow the lines of the proof given for $r=1$ by Dupont and
Vigu\'e-Poirrier in \cite{DV}.   Write $\land V = \land
V^{\mbox{\scriptsize even}} \otimes \land V^{\mbox{\scriptsize
odd}}$, and denote by $(x_i)_{i\in I}$ a graded basis of $
V^{\mbox{\scriptsize even}} \oplus V^{\mbox{\scriptsize odd}}$. We
denote by $\frac{\partial }{\partial x_i}$ the derivation of degree
$-\vert x_i\vert$ defined by
$$\frac{\partial}{\partial x_i} (x_i) = 1 \, \, \mbox{and } \frac{\partial}{\partial x_i} (x_j) =
0\,, i\neq j\,.$$

If $v \in V_p^q$, we denote $\ell (v) = p+q$. This is a new
gradation, and for any element $P$ of $\land V$, we have
$$\ell(P)\, P = \sum_i \, \ell(x_i) \, x_i\, \frac{\partial}{\partial x_i} (P)\,.$$

The lower gradation on $V$ extends to $\overline{V}$. If $v\in
V_p^q$, then $s(v)\in \overline{V}^{q-r}_p$. The differential
$\overline{D}$ is compatible with this double gradation,
$$\overline{D} : (\land V\otimes  \overline{V})^q_p \to
(\land V\otimes \overline{V})_{p-1}^{q+1}\,.$$

 Write $P= \overline{D}x $, $P_i =  \overline{D}x_i$ and $\omega = \sum \overline{x_i}\,
a_i$ with $x_i\in V$, $a_i\in \land^+V$. Then
$$0 =  \overline{D}\overline{v} + \sum  \overline{D}
(\overline{x_i} a_i) = (-1)^r \left(s(P) + \sum_i s(P_i)\,
a_i\right) + \sum_i (-1)^{\vert \overline{x_i}\vert}
\overline{x_i}\cdot \overline{D}(a_i)$$
$$= (-1)^r \left( \sum_i \overline{x_i}\, \frac{\partial P}{\partial x_i} +
 \sum_{ij} \overline{x_i} \,\frac{\partial P_j}{\partial x_i}\, a_j\right)+
 \sum_i (-1)^{\vert \overline{x_i}\vert} \overline{x_i}\cdot  \overline{D}(a_i)\,.$$
Therefore
$$\frac{\partial P}{\partial x_i} = -(-1)^{\vert x_i\vert} \overline{D}a_i -
 \sum_j \frac{\partial P_j}{\partial x_i} \, a_j\,,$$
and
$$\ell (P) P = \sum_i \ell (x_i) x_i \frac{\partial P}{\partial x_i}
 = -\left( \sum_{ij} \ell (x_i) x_i\, \frac{\partial P_j}{\partial x_i}\, a_j
 + \sum_i (-1)^{\vert x_i\vert} \ell (x_i) x_i \,  \overline{D}a_i\right)$$
$$= - \left(\sum_i \ell (P_i) P_i\, a_i + \sum_i (-1)^{x_i} \ell (x_i) x_i \,  \overline{D}a_i\right) =
-  \overline{D}\left( \sum_i \ell (x_i) x_i\, a_i\right)\,.$$ This
implies that $$v+ \sum_i \frac{\ell (x_i)}{\ell (x)} \, x_i \, a_i$$
is a cocycle. In particular, $v\in V_0$ and is a cocycle. This ends
the proof of   theorem~2.

 \end{document}